\documentclass[aap,seceqn,nameyear,MSNbibl,dvips]{arximspdf}
\usepackage{mathrsfs}

\doi{10.1214/10-AAP748}
\volume{21}
\issue{5}
\pubyear{2011}
\firstpage{1994}
\lastpage{2015}

\makeatletter

\def\bptnote#1{}

\newtheorem{theorem}{Theorem}
\newtheorem{proposition}{Proposition}
\newtheorem{lemma}{Lemma}
\newproclaim{remark}{Remark}

\newtheorem{corollary}{Corollary}

\newcommand{\X}{{\mathbf X}}

\newcommand{\ep}{\varepsilon}
\newcommand{\ga}{\alpha}

\newcommand{\gd}{\delta}
\newcommand{\gS}{\Sigma}
\newcommand{\bga}{\bolds\alpha}
\newcommand{\rtr}{\operatorname{tr}}
\newcommand{\rE}{\mathrm{E}}
\newcommand{\fn}{\frac{1}{n}}
\newcommand{\skln}{\sum_{k=1}^n}
\newcommand{\sjln}{\sum_{j=1}^n}

\newcommand{\bU}{\mathbf{U}}
\newcommand{\bs}{\mathbf{s}}
\newcommand{\bu}{\mathbf{u}}
\newcommand{\bz}{\mathbf{z}}
\newcommand{\by}{\mathbf{y}}
\newcommand{\bx}{\mathbf{x}}
\newcommand{\bV}{\mathbf{V}}
\newcommand{\bt}{\mathbf{t}}
\newcommand{\rP}{\mathrm{P}}
\newcommand{\bbA}{\mathbf{A}}
\newcommand{\bbB}{\mathbf{B}}
\newcommand{\bbe}{\mathbf{e}}
\newcommand{\bbI}{\mathbf{I}}
\newcommand{\bbj}{\mathbf{j}}
\newcommand{\bbM}{\mathbf{M}}
\newcommand{\bbS}{\mathbf{S}}
\newcommand{\bbT}{\mathbf{T}}
\newcommand{\bbU}{\mathbf{U}}
\newcommand{\bbv}{\mathbf{v}}
\newcommand{\bbX}{\mathbf{X}}
\newcommand{\bbY}{\mathbf{Y}}
\newcommand{\bbz}{\mathbf{z}}
\newcommand{\cC}{\mathcal{C}}
\newcommand{\cF}{\mathcal{F}}
\newcommand{\gma}{\gamma}

\makeatother

\begin{document}
\begin{frontmatter}

\title{Asymptotic properties of eigenmatrices of a large sample
covariance matrix}
\runtitle{Eigenmatrices of a sample covariance matrix}

\begin{aug}
\author[A]{\fnms{Z. D.} \snm{Bai}\corref{}\thanksref{t1}\ead[label=e1]{baizd@nenu.edu.cn}},
\author[B]{\fnms{H. X.} \snm{Liu}\ead[label=e2]{huixialiu@dbs.com.sg}} and
\author[C]{\fnms{W. K.} \snm{Wong}\thanksref{t2}\ead[label=e3]{ecswwk@nus.edu.sg}}
\runauthor{Z. D. Bai, H. X. Liu and W. K. Wong}
\affiliation{Northeast Normal University, National University of
Singapore\break and Hong Kong Baptist University}
\address[A]{Z. D. Bai\\
KLASMOE and School of Mathematics\\
\quad and Statistics\\
Northeast Normal University\\
Changchun 130024\\
China\\
\printead{e1}}
\address[B]{H. X. Liu\\
Department of Statistics\\
\quad and Applied Probability\\
National University of Singapore\\
Science dr. 2 117546\\
Singapore\\
\printead{e2}}
\address[C]{W. K. Wong\\
Department of Economics\\
Hong Kong Baptist University\\
Hong Kong\\
China\\
\printead{e3}}
\end{aug}

\thankstext{t1}{Supported by the NSFC
10871036 and NUS Grant R-155-000-079-112.}

\thankstext{t2}{Supported by grants from Hong Kong Baptist University.}

\received{\smonth{10} \syear{2009}}
\revised{\smonth{2} \syear{2010}}

%
\begin{abstract}
Let $S_n=\frac1 n X_nX_n^*$ where $X_n=\{X_{ij} \}$ is a $p\times
n$ matrix with i.i.d. complex standardized entries having finite fourth
moments. Let $Y_n(\bt_1,\bt_2,\sigma)=
\sqrt{p}({\bx}_n(\bt_1)^*(S_n+\sigma
I)^{-1}{\bx}_n(\bt_2)-{\bx}_n(\bt_1)^*{\bx}_n(\bt_2) m_n(\sigma))$ in
which $\sigma>0$ and $m_n(\sigma)=\int\frac{dF_{y_n}(x)}{x+\sigma
}$ where
$F_{y_n}(x)$ is the Mar\v{c}enko--Pastur law with parameter
$y_n=p/n$; which converges to a positive constant as $n\to\infty$,
and ${\bx}_n(\bt_1)$ and ${\bx}_n(\bt_2)$ are unit
vectors in ${\Bbb
C}^p$, having indices $\bt_1$ and $\bt_2$, ranging in a compact
subset of a finite-dimensional Euclidean space. In this paper, we
prove that the sequence $Y_n(\bt_1,\bt_2,\sigma)$ converges
weakly to
a $(2m+1)$-dimensional Gaussian process. This result provides
further evidence in support of the conjecture that the distribution
of the eigenmatrix of $S_n$ is asymptotically close to that of a
Haar-distributed unitary matrix.
\end{abstract}

%

\end{frontmatter}

\section{Introduction}
\label{intro}

Suppose that $\{x_{jk}, j,k=1,2,\ldots\}$ is a double array of
complex random variables that are independent and identically
distributed (i.i.d.) with mean zero and variance $1$. Let
${\bx}_j=(x_{1j},\ldots, x_{pj})'$ and $\X=({\bx}_{1},\ldots,
{\bx}_{n})$, we define
%
%
\begin{equation}\label{sn}
\bbS_n=\frac1 n \sum_{k=1}^n{\bx}_k{\bx}_k^*=\frac1 n
\X\X^*,
\end{equation}
where ${\bx}_k^*$ and $\X^*$ are the transposes of the complex
conjugates of ${\bx}_k$ and $\X$, respectively. The matrix $\bbS_n$
defined in (\ref{sn}) can be viewed as the sample covariance matrix
of a $p$-dimensional random sample with size $n$. When the dimension
$p$ is fixed and the sample size $n$ is large, the spectral behavior
of $\bbS_n$ has been extensively investigated in the literature due
to its importance in multivariate statistical inference [see, e.g.,
Anderson (\citeyear{And51}, \citeyear{And89})]. However, when the
dimension $p$ is
proportional to the sample size $n$ in the limit; that is, $\frac
p n\rightarrow y>0$ as $n\rightarrow\infty$, the classical
asymptotic theory will induce serious inaccuracy. This phenomenon
can be easily explained from the viewpoint of random matrix theory
(RMT).

Before introducing our advancement of the theory, we will first give
a brief review of some well-known properties of $\bbS_n$ in RMT. We
define the empirical spectral distribution (ESD) of $\bbS_n$ by
\[
F^{\bbS_n}(x)=\frac{1}{p}\sum_{j=1}^pI(\lambda_j\leq x),
\]
where $\lambda_j$'s are eigenvalues of $\bbS_n$.
First, it has long been known that $F^{\bbS_n}(x)$ converges almost
surely to the standard Mar\v{c}enko--Pastur law [MPL; see, e.g.,
\citet{MP67}, \citet{Wac78} and \citet{Yin86}]
$F_y(x)$, which has a density
function $(2\pi xy)^{-1}\sqrt{(b-x)(x-a)}$, supported on
$[a,b]=[(1-\sqrt{y})^2, (1+\sqrt{y})^2]$. For the case $y>1$,
$F_y(x)$ has a point mass $1-1/y$ at 0. If its fourth moment is
finite, as $n\rightarrow\infty$, the largest eigenvalue of $\bbS_n$
converges to $b$ while the smallest eigenvalue (when $y \le1$) or
the $(p-n+1)$st smallest eigenvalue (when $y>1$) converges to
$a$ [see \citet{Bai99} for a review]. The central limit theorem
(CLT) for
linear spectral statistics (LSS) of $\bbS_n$ has been
established in \citet{BaiSil04}.

While results on the eigenvalues of $\bbS_n$ are abundant in the
literature, not much work has been done on the behavior of the
eigenvectors of $\bbS_n$. It has been conjectured that the
eigenmatrix; that is, the matrix of orthonormal eigenvectors of
$\bbS_n$, is asymptotically Haar-distributed. This conjecture has
yet to be formally proven due to the difficulty of describing the
``asymptotically Haar-distributed'' properties when the dimension
$p$ increases to infinity. \citet{Sil81} was the first one to
create an approach to characterize the eigenvector properties. We
describe his approach as follows: denoting the spectral
decomposition of $\bbS_n$ by $\bbU_n^*{\bolds\Lambda}
\bbU_n$, if ${\bx}_{ij}$ is normally distributed, $\bbU_n$ has
a Haar
measure on the orthogonal matrices and is independent of the
eigenvalues in ${\bolds\Lambda}$. For any unit vector
${\bx}_n\in
{\Bbb C}^p$, the vector $\by_n=(y_1,\ldots,y_p)=\bbU_n{\bx}_n$ performs
like a uniform distribution over the unit sphere in ${\Bbb C}^p$. As
such, for $t\in[0,1]$, a stochastic process
\[
X_n(t)=\sqrt{p/2}\Biggl(\sum_{i=1}^{[pt]}(y_i^2-1/p)\Biggr)
\]
is defined. If $\bz=(z_1,\ldots,z_p)'\sim N(0,I_p)$, then $\by_n$ has
the same
distribution as $\bz/\|\bz\|$ and $X_n(t)$ is identically distributed with
\[
\tilde X_n(t)=\sqrt{p/2}\|\bz\|^{-2}\Biggl(\sum_{i=1}^{[pt]}(z_i^2-\|\bz
\|
^2/p)\Biggr).
\]
Applying Donsker's theorem [\citet{Don51}], $X_n(t)$ tends to a
standard Brownian bridge.

For any general large sample covariance, it is important to examine
the behavior of the $X_n(t)$ process. Silverstein (\citeyear{Sil81},
\citeyear{Sil84}, \citeyear{Sil89})
prove that the integral of polynomial functions
with respect to $X_n(t)$ will tend to a normal distribution. To
overcome the difficulty of tightness, \citet{Sil90} takes
${\bx}_n=(\pm1,\ldots,\pm1)/\sqrt{p}$ so that the process $X_n(t)$
will tend to the standard Brownian bridge instead. In addition,
\citet{BaiMiaPan07} investigate the process $X_n(t)$, defined for
$\bbT_n^{1/2}\bbS_n\bbT_n^{1/2}$ with $(\bbT_n^{1/2})^2=\bbT_n$, a
nonnegative positive definite matrix.

However, so far, the process $X_n(t)$ is assumed to be generated
only by one unit vector ${\bx}_n$ in ${\Bbb C}^p$. This imposes
restrictions on many practical situations. For example, in the
derivation of the limiting properties of the bootstrap corrected
Markowitz portfolio estimates, we need to consider two unparallel
vectors simultaneously [see \citet{BaiLiuWon09} and
Markowitz (\citeyear{M52}, \citeyear{Mar59}, \citeyear{M91})]. In this
paper, we will go beyond the boundaries
of their studies to investigate the asymptotics of the eigenmatrix
for any general large sample covariance matrix $\bbS_n$ when ${\bx}_n$
runs over a subset of the $p$-dimensional unit sphere in which
${\Bbb C}_1^p=\{{\bx}_n\dvtx\|{\bx}_n\|=1, {\bx}_n\in
{\Bbb C}^p\}$.

We describe the approach we introduced in this paper as follows:
if $\bV_n$ is Haar-distributed, for any pair
of $p$-vectors $\bx$ and $\by$ satisfying $\bx\perp\by$, $(\bV_n\bx, \bV
_n\by)$ possesses the same joint
distribution as
%
%
\begin{equation} \label{y1}
\pmatrix{\bz_1&\bz_2\cr
}
\pmatrix{ \bz_1^*\bz_1&\bz_1^*\bz_2\cr\bz_2^*\bz_1&\bz_2^* \bz_2\cr
}
^{-1/2} ,
\end{equation}
where $\bz_1$ and $\bz_2$ are two
independent $p$-vectors whose components are i.i.d. standard normal
variables. As $n$ tends to infinity, we have
%
%
\begin{equation} \label{y2}
\frac1 p
\pmatrix{\bz_1^*\bz_1&\bz_1^*\bz_2\cr\bz_2^*\bz_1&\bz_2^*\bz_2\cr
}
\longrightarrow
I_2.
\end{equation}
Therefore, any group of functionals defined by
these two random vectors should be asymptotically independent of
each other. We shall adopt this setup to explore
the conjecture that $\bU_n$ is asymptotically Haar-distributed.

We consider ${\bx}$ and $\by$ to be two $p$-vectors with an angle
$\theta$. Thereafter, we find two orthonormal vectors $\bga_1$ and
$\bga_2$ such that
\[
{\bx}=\|{\bx}\|\bga_1 \quad\mbox{and}\quad
\by=\|\by\|(\bga_1\cos\theta+\bga_2\sin\theta).
\]
By (\ref{y1}) and (\ref{y2}), we have
%
%
\begin{equation} \label{y3}
\bV_n\bx\sim
p^{-1/2}\|{\bx}\|\bz_1 \quad\mbox{and}\quad \bV_n\by\sim
p^{-1/2}\|\by\|(\bz_1\cos\theta+\bz_2\sin\theta).
\end{equation}
Let
$\sigma>0$ be a positive constant, we now consider the following three
quantities:
%
%
\begin{equation}\label{threequan}
{\bx}^*(\bbS_n+\sigma I)^{-1}{\bx},\qquad
{\bx}^*(\bbS_n+\sigma
I)^{-1}\by\quad\mbox{and}\quad\by^*(\bbS_n+\sigma I)^{-1}\by.
\end{equation}

We hypothesize that if $\bU_n$ is asymptotically Haar-distributed and
is asymptotically independent of ${\bolds\Lambda}$, then the
above three quantities should be asymptotically equivalent to
%
%
\begin{eqnarray} \label{cltq}
&p^{-1}\|{\bx}\|^2\bz_1^*({\bolds\Lambda}+\sigma
I)^{-1}\bz_1,&\nonumber\\
&p^{-1}\|{\bx}\|\|\by\|\bz_1^*({\bolds\Lambda}+\sigma
I)^{-1}(\bz_1\cos\theta+\bz_2\sin\theta)\quad\mbox{and}& \\
&p^{-1}\|\by\|^2(\cos\theta
\bz_1+\sin\theta\bz_2)^*({\bolds\Lambda}+\sigma
I)^{-1}(\bz_1\cos\theta+\bz_2\sin\theta),&\nonumber
\end{eqnarray}
respectively. We then proceed to investigate the stochastic
processes related to these functionals. By using the Stieltjes
transform of the sample covariance matrix, we have
\[
p^{-1}\bz_1^*({\bolds\Lambda}+\sigma I)^{-1}\bz_1\to
m(\sigma)=-\frac{1+\sigma-y -\sqrt{(1+y+\sigma)^2-4y}}{2y\sigma}
\qquad\mbox{a.s.},
\]
where $m(\sigma)$ is a solution to the quadratic equation
%
%
\begin{equation} \label{meq}
m(1+\sigma
-y+y\sigma m)-1=0.
\end{equation}
Here, the selection of $m(\sigma)$ is
due to the fact that $m(\sigma)\to0$ as
$\sigma\to\infty$. By using the same argument, we conclude that
\[
p^{-1}(\cos\theta\bz_1+\sin\theta\bz_2)^*({\bolds\Lambda}+\sigma
I)^{-1}(\bz_1\cos\theta+\bz_2\sin\theta)\longrightarrow m(\sigma)\qquad
\mbox{a.s.}
\]
Applying the results in \citet{BaiMiaPan07}, it can be easily
shown that, for the complex case,
%
%
\begin{equation} \label{clt1}
p^{-1/2}[\bz_1^*({\bolds\Lambda}+\sigma I)^{-1}\bz
_1-pm_n(\sigma)]
\longrightarrow N(0,W),
\end{equation}
and for the real case, the
limiting variance is $2W$, where $W=W(\sigma,\sigma)$, $m_n(\sigma)$ is
$m(\sigma)$ with $y$ replaced by $y_n$ such that
\begin{eqnarray*}
m_n(\sigma)&=&-\frac{1+\sigma-y_n -\sqrt{(1+y_n+\sigma
)^2-4y_n}}{2y_n\sigma},\\
y_n&=&p/n
\end{eqnarray*}
and
\[
W(\sigma_1,\sigma_2)=\frac{m(\sigma_1)m(\sigma_2)}{1-y(1-\sigma
_1ym(\sigma_1))(1-\sigma_2m(\sigma_2))}.
\]
Here, the definitions of ``real case'' and ``complex
case'' are given in Theorem \ref{thm1} as stated in the next
section. By the same argument, one could obtain a similar result
such that
%
%
\begin{equation} \label{clt2}\qquad
p^{-1/2}[(\bz_1\cos\theta
+\bz_2\sin\theta)^*({\bolds\Lambda}+\sigma
I)^{-1}(\bz_1\cos\theta+\bz_2\sin\theta)-pm_n(\sigma)].
\end{equation}

We normalize the second term in (\ref{cltq}) and, thereafter, derive
the CLT for the joint distribution of all three terms stated in
(\ref{cltq}) after normalization. More notably, we establish some
limiting behaviors of the processes defined by these normalized
quantities.

\section{Main results}
\label{sec:1}

Let $\Bbb S=\Bbb S_p$ be a subset of the unit $p$-sphere ${\Bbb
C}_1^p$ indexed by an $m$-dimensional hyper-cube $T=[0,2\pi]^m$.
For any $m$ arbitrarily chosen orthogonal unit
$p$-vectors ${\bx}_1,\ldots,{\bx}_{m+1}\in{\Bbb C}_1^p$,
we define
%
%
\begin{eqnarray} \label{s}\qquad
{\Bbb S} &=& \{{\bx}_n(\bt)={\bx}_1\cos
t_1+{\bx}_2\sin t_1\cos
t_2+\cdots+{\bx}_{m}\sin t_1\cdots\sin t_{m-1}\cos t_{m}\nonumber\\
&&\hspace*{135.5pt}{} +{\bx}_{m+1}\sin t_1\cdots\sin t_{m-1}\sin
t_{m},
\bt\in T \}.
\end{eqnarray}
If $\Bbb S$ is chosen in the form of (\ref{s}), then the inner
product $\bx_n(\bt_1)^*\bx_n(\bt_2)$ is a function of
$\bt_1$
and $\bt_2$ only (i.e., independent of $n$). Also, the norm
of the difference (we call it norm difference in this paper)
$\|\bx_n(\bt_1)-\bx_n(\bt_2)\|$ satisfies the Lipschitz
condition. If the time index set is chosen arbitrarily, we could
assume that the angle, $\vartheta_n(\bt_1,\bt_2)$, between
${\bx}_n(\bt_1)$ and $\bx_n(\bt_2)$ tends to a
function of $\bt_1$
and $\bt_2$ whose norm difference satisfies the Lipschitz
condition.

Thereafter, we define a stochastic process $\bbY_n(\bu,\sigma)$
mapping from the time index set $T\times T\times I$ to ${\Bbb S}$
with $I=[\sigma_{10},\sigma_{20}]$ ($0<\sigma_{10}<\sigma_{20}$)
such that
\[
Y_n(\bu,\sigma)= \sqrt{p}\bigl({\bx}_n(\bt_1)^*(\bbS
_n+\sigma
I)^{-1}{\bx}_n(\bt_2)-\bx_n(\bt_1)^*\bx_n(\bt_2) m_n(\sigma)\bigr),
\]
where $(\bu,\sigma)=(\bt_1,\bt_2,\sigma)\in T\times
T\times I$.
\begin{remark} \label{remsphere} If the sample covariance matrix $\bbS
_n$ is real,
the vectors $\bx_n$ and $\by_n$ will be real,
and thus, the set $\Bbb S$ has to be defined as a subset of unit
sphere $\Bbb R^p_1=\{\bx\in\Bbb R^p, \|\bx\|=1\}$. The time index
can be similarly described for the complex case. In what follows, we
shall implicitly use the convention for the real case.
\end{remark}

We have the following theorem.
\begin{theorem}\label{thm1}
Assume that the entries of $\bbX$ are i.i.d. with mean 0, variance 1,
and finite fourth moments. If the variables are complex, we further
assume $\rE X_{11}^2=0$ and $\rE|X_{11}|^4=2$, and refer to this
case as the \textup{complex case}. If the variables are real, we assume
$\rE X_{11}^4=3$ and refer to it as the \textup{real case}. Then, as
$n\rightarrow\infty$, the process $Y_n(\bt_1,\bt_2,\sigma)$ converges
weakly to a multivariate Gaussian process $Y(\bt_1,\bt_2,\sigma)$
with mean zero and variance--covariance function $\rE
Y(\bt_1,\bt_2,\sigma_1)Y(\bt_3,\bt_4,\sigma_2)$ satisfying
\[
\rE
Y(\bt_1,\bt_2,\sigma_1)Y(\bt_3,\bt_4,\sigma
_2)=\vartheta(\bt_1,\bt_4)\vartheta(\bt_3,\bt_2)
W(\sigma_1,\sigma_2)
\]
for the complex case and satisfying
\[
\rE
Y(\bt_1,\bt_2,\sigma_1)Y(\bt_3,\bt_4,\sigma_2) =
\bigl(\vartheta(\bt_1,\bt_4)\vartheta(\bt_3,\bt_2)
+\vartheta(\bt_1,\bt_3)\vartheta(\bt_4,\bt_2)
\bigr)W(\sigma_1,\sigma_2)
\]
for the real case where
\[
W(\sigma_1,\sigma_2)=\frac{ym(\sigma_1)m(\sigma_2)}{1-y(1-\sigma
_1m(\sigma_1))(1-\sigma_2m(\sigma_2))}
\]
and
\[
\vartheta(\bt,\bs)=\lim\bx_n^*(\bt)\bx_n(\bs).
\]
\end{theorem}

We will provide the proof of this theorem in the next section. We
note that \citet{BaiMiaPan07} have proved that
\[
\sqrt{p}[{\bx}_n(\bt_1)^*(\bbS_n+\sigma
I)^{-1}{\bx}_n(\bt_1)-m_n(\sigma)] \longrightarrow N(0,W)
\]
for the complex case and proved that the
asymptotic variance is $2W$ for the real case.

More generally, if ${\bx}$ and $\by$ are two orthonormal vectors,
applying Theorem \ref{thm1}, we obtain the limiting distribution of
the three quantities stated in (\ref{threequan})
with normalization such that
%
%
\begin{equation} \label{threequanorm}\qquad
\sqrt{p}
\pmatrix{{\bx}^*(\bbS_n+\sigma I)^{-1}{\bx}-m_n(\sigma)\cr
{\bx}^*(\bbS_n+\sigma I)^{-1}\by\cr\by^*(\bbS_n+\sigma
I)^{-1}\by-m_n(\sigma)\cr
}
\longrightarrow N
\left(
\pmatrix{0\cr0\cr0
}
,
\pmatrix{W&0&0\cr
0&W&0\cr0&0&W\cr
}
\right)
\end{equation}
for
the complex case while the asymptotic covariance matrix is
\[
\pmatrix{2W&0&0\cr0&W&0\cr0&0&2W\cr
}
\]
for the real case.
\begin{remark}
This theorem shows that the three quantities stated in
(\ref{threequan}) are asymptotically independent of one another.
It provides a stronger support to the conjecture that $\bU_n$ is
asymptotically Haar-distributed than those established in the
previous literature.
\end{remark}

In many practical applications, such as wireless communications and
electrical engineering [see, e.g., \citet{EvaTse00}], we are
interested in extending the process $Y_n(\bu,\sigma)$
defined on a region $T\times T\times D$ where $D$ is a compact
subset of the complex plane and is disjoint with the interval
$[a,b]$, the support of the MPL. We can define a complex measure by
putting complex mass
$\bx^*(\bt_1)\bbU_n^*\bbe_j\bbe_j'\bbU_n\by(\bt_2)$
at $\lambda_j$,
the $j$th eigenvalue of $\bbS_n$, where $\bbe_j$ is the
$p$-vector with 1 in its $j$th entry and 0 otherwise. In this
situation, the Stieltjes transform of this complex measure is
\[
s_n(z)=\bx^*(\bbS_n-z \bbI)^{-1}\by,
\]
where $z=\mu+iv$ with $v\neq0$. When considering the CLT of LSS
associated with the complex measure defined above, we need to
examine the limiting properties of the Stieltjes transforms, which
lead to the extension of the process $Y_n(\bu,\sigma)$ to
$Y_n(\bu,-z)$, where $z$ is an index number in $D$.

If ${\bx}^*\by$ is a constant (or has a limit, we
still denote it as ${\bx}^*\by$ for simplicity), it follows from Lemma
\ref{lema1} that
\[
\bx^*(\bbS_n-z \bbI)^{-1}\by\longrightarrow{\bx}^*\by s(z),
\]
where
\begin{eqnarray*}s(z)&=&
\cases{\displaystyle\frac{1-z-y+\sqrt{(1-z+y)^2-4y}}{2yz},
&\quad when $\Im(z)>0$,\vspace*{2pt}\cr
\bar s(\bar z), &\quad when $\Im(z)<0$,}
\\
&=&\frac{1-z-y+\operatorname{sgn}(\Im(z))\sqrt{(1-z+y)^2-4y}}{2yz}
\qquad\mbox{if
} \Im(z)\ne0,
\end{eqnarray*}
is the Stieltjes transform of MPL, in which, by
convention, the square root $\sqrt z$ takes the one with the
positive imaginary part. When $z\ne0$ is real, $s(z)$ is defined as
the limit from the upper complex plane. By definition,
$m(\sigma)=s(-\sigma+i0)=\lim_{v\downarrow0}s(-\sigma+iv)$. In calculating
the limit, we follow the conventional sign of the square root of a
complex number that the real part of $\sqrt{(-\sigma+iv-1-y)^2-4y}$
should have the opposite sign of $v$, and thus
\[
m(\sigma)=-\frac{1+\sigma-y-\sqrt{(1+y+\sigma)^2-4y}}{2y\sigma}.
\]
Now, we are ready to extend the process $Y_n(\bu,\sigma)$ to
\[
Y_n(\bu,z)=\sqrt{p}[\bx_n^*(\bt_1)(\bbS_n-z
\bbI)^{-1}\bx_n(\bt_2)-\bx_n^*(\bt_1)\bx_n(\bt_2)s(z,y_n)],
\]
where $s(z,y_n)$ is the Stieltjes transform of the LSD of $\bbS_n$
in which $y$ is replaced by $y_n$. Here, $z=u+iv$
with $v>0$ or $v<0$. Thereby, we obtain the following theorem.
\begin{theorem}\label{thme2}
Under the conditions of Theorem \ref{thm1}, the process
$Y_n(\bu,z)$ tends to a multivariate Gaussian process $Y(\bu,z)$
with mean 0 and covariance function $\rE(Y(\bu,z_1)Y(\bu,z_2))$
satisfying
%
%
\begin{equation}\rE(Y(\bu,z_1)Y(\bu,z_2))= \vartheta(\bt_1,\bt_4)
\vartheta(\bt_3,\bt_2)W(z_1,z_2)
\end{equation}
for the complex case and
satisfying
%
%
\begin{equation} \label{varreal}\qquad
\rE(Y(\bu,z_1)Y(\bu,z_2))=
\bigl(\vartheta(\bt_1,\bt_4)
\vartheta(\bt_3,\bt_2)+\vartheta(\bt_1,\bt_3)
\vartheta(\bt_4,\bt_2)\bigr)W(z_1,z_2)
\end{equation}
for the real
case where
\[
W(z_1,z_2)=\frac{ys(z_1)s(z_2)}{1-y(1+z_1s(z_1))(1+z_2s(z_2))}.
\]
\end{theorem}

Theorem \ref{thme2} follows from Theorem \ref{thm1} and Vitali lemma
[see Lemma 2.3 of \citet{BaiSil04}] since both $Y(\bu,z)$
and $Y_n(\bu,z)$ are analytic functions when $z$ is away from
$[a,b]$, the support of MPL.

Suppose that $f(x)$ is analytic on an open region containing the
interval $[a,b]$. We construct an LSS with respect to the complex
measure as defined earlier; that is,
\[
\sum_{j=1}^p
f(\lambda_j)\bx^*(\bt_1)\bbU_n^*\bbe_j\bbe_j'\bbU_n\by(\bt_2) .
\]
We then consider the normalized quantity
%
%
\begin{eqnarray}\label{lls1}
X_n(f)&=&\sqrt{p}\Biggl(\sum_{j=1}^p
f(\lambda_j)\bx^*(\bt_1)\bbU_n^*\bbe_j\bbe_j'\bbU_n\by(\bt_2)\nonumber\\
[-8pt]\\[-8pt]
&&\hspace*{27.14pt}{}-{\bx}^*(\bt_1)\by(\bt_2)\int
f(x)\,dF_{y_n}(x)\Biggr),\nonumber
\end{eqnarray}
where $F_y$ is the
standardized MPL. By applying the Cauchy formula
\[
f(x)=\frac1{2\pi i}\oint_{\cC}\frac{f(z)}{z-x}\,dz,
\]
where $\cC$ is a contour enclosing $x$, we obtain
%
%
\begin{eqnarray}\label{lls2}
X_n(f,\bu)&=&-\frac{\sqrt{p}}{2\pi
i}\biggl(\oint_{\cC}[\bx_n^*(\bt_1)(\bbS_n-z\bbI)^{-1}\by
(\bt_2)\nonumber\\[-8pt]\\[-8pt]
&&\hspace*{66.23pt}{}
-{\bx}_n^*(\bt_1)\by(\bt_2)s_n(z)]f(z)\,dz\biggr),\nonumber
\end{eqnarray}
where
$\cC$ is a contour enclosing the interval $[a,b]$,
$\bu=(\bt_1,\bt_2)$, and
\[
s_n(z)=\frac{1-z-y_n+\operatorname{sgn}(\Im(z))\sqrt{(1-z+y_n)^2-4y_n}}{2y_nz}.
\]
Thereafter, we obtain the following two corollaries.
\begin{corollary}\label{coro1}
Under the conditions of Theorem \ref{thm1}, for any $k$ functions
$f_1,\ldots,f_k$ analytic on an open region containing the interval
$[a,b]$, the $k$-dimensional process
\[
(X_n(f_1,\bu_1),\ldots,X_n(f_k,\bu_k))
\]
tends to the $k$-dimensional stochastic multivariate Gaussian
process with mean zero and covariance function satisfying
\[
\rE(X(f,\bu)X(g,\bbv))=-\frac\theta
{4\pi^2}\oint_{c_1}\oint_{c_2}W(z_1,z_2)f(z_1)g(z_2)\,dz_1\,dz_2,
\]
where $\theta=\vartheta(\bt_1,\bt_4) \vartheta(\bt_3,\bt_2)$
for the complex case and $=\vartheta(\bt_1,\bt_4)
\vartheta(\bt_3,\bt_2)+\vartheta'(\bt_1,\bt_3)
\overline{\vartheta}{}'(\bt_4,\bt_2)$ for the real case. Here,
$\cC_1$ and $\cC_2$ are two disjoint contours that enclose the
interval $[a,b]$ such that the functions $f_1,\ldots,f_k$ are
analytic inside and on them.
\end{corollary}
\begin{corollary}\label{coro2}
The covariance function in Corollary \ref{coro1} can also be
written as
\begin{eqnarray*}
\rE(X(f,\bu)X(g,\bbv))
&=& \theta
\biggl(\int_{a}^bf(x)g(x)\,dF_y(x)\\
&&\hspace*{9.5pt}{}-\int_{a}^bf(x)\,dF_y(x)\int
_{a}^bg(x)\,dF_y(x)\biggr),
\end{eqnarray*}
where $\theta$ has been defined in Corollary \ref{coro1}.
\end{corollary}

\section{\texorpdfstring{The proof of Theorem \protect\ref{thm1}}{The proof of Theorem 1}}
\label{sec:2}

To prove Theorem \ref{thm1}, by Lemma \ref{lem33}, it is
sufficient to show that $Y_n(\bu,\sigma)-\rE Y_n(\bu,\sigma)$
tends to
the limit process $Y(\bu,\sigma)$. We will first prove
the property of the finite-dimensional convergence in Section
\ref{secfin} before proving the tightness property in
Section \ref{sectight}. Throughout the paper, the limit is
taken as $n\rightarrow\infty$.

\subsection{Finite-dimensional convergence}\label{secfin}

Under the assumption of a finite fourth moment, we follow \citet
{BaiMiaPan07} to truncate the random variables $X_{ij}$ at
$\ep_n\sqrt[4]{n}$ for all $i$ and $j$ in which $\ep_n\to0$ before
renormalizing the random variables to have mean 0 and variance 1.
Therefore, it is reasonable to impose an additional assumption that
$ |X_{ij}|\le\ep_n\sqrt[4]{n}$ for all $i$ and $j$.

Suppose $\bs_j$ denotes the $j$th column of
$\frac{1}{\sqrt{n}}X_n$. Let $\bbA(\sigma)=\bbS_n+\sigma\bbI$ and
$\bbA_j(\sigma)=\bbA(\sigma)-\bs_j\bs_j^*$. Let ${\bx}_n$ and $\by_n$
be any
two vectors in $\Bbb C_1^p$. We define
\begin{eqnarray*}
\xi_j(\sigma)&=&\bs_j^{*}\bbA_j^{-1}(\sigma)\bs_j-\frac
{1}{n}\rtr
\bbA_j^{-1}(\sigma), \\
\gamma_j
&=&\bs_j^*\bbA_j^{-1}\by_n\bx_n^*\bbA_j^{-1}(\sigma)\bs_j-\frac{1}{n}\bx
_n^*\bbA_j^{-1}(\sigma)\bbA_j^{-1}(\sigma
)\by_n
, \\
\beta_j(\sigma)&=&\frac{1}{1+\bs_j^*\bbA_j^{-1}(\sigma)\bs_j},\\
b_j(\sigma)&=&\frac{1}{1+n^{-1}\rtr\bbA_j^{-1}(\sigma)}
\end{eqnarray*}
and
\[
\bar b
=\frac1{1+n^{-1}\rE\rtr\bbA^{-1}(\sigma)}.
\]
We also define the
$\sigma$-field ${\cF}_j=\sigma(\bs_1,\ldots,\bs_j)$. We denote
by $\rE_j(\cdot)$ the conditional expectation when $\cF_j$ is given. By
convention, $\rE_0$ denotes the unconditional expectation.

Using the martingale decomposition, we have
%
%
\begin{eqnarray} \label{martina}
\bbA^{-1}(\sigma)-\rE
\bbA^{-1}(\sigma)&=&\sjln
(\rE_j-\rE_{j-1})[\bbA^{-1}(\sigma)-\bbA_k^{-1}(\sigma)]\nonumber\\
[-8pt]\\[-8pt]
&=&\sjln
(\rE_j-\rE_{j-1})\beta_j\bbA_j^{-1}(\sigma)\bs_j\bs_j^*
\bbA_j^{-1}(\sigma).\nonumber
\end{eqnarray}
Therefore,
\begin{eqnarray*}
Y_n(\bu,\sigma)&=&\sqrt{p}\sjln
(\rE_j-\rE_{j-1})\bx(\bt_1)^*[\bbA^{-1}(\sigma)-\bbA
_k^{-1}(\sigma)]\bx(\bt_2)\\
&=&\sqrt{p}\sjln
(\rE_j-\rE_{j-1})\beta_j\bx(\bt_1)^*\bbA_j^{-1}(\sigma
)\bs_j\bs_j^*
\bbA_j^{-1}(\sigma)\bx(\bt_2).
\end{eqnarray*}

Consider the $K$-dimensional distribution of $\{Y_n(\bu_1,\sigma
_1),\ldots, Y_n(\bu_K,\sigma_K)\}$\break where
$(\bu_i,\sigma_i)=(\bt_{i1},\bt_{i2}, \sigma_i)\in
T\times T\times
I$. Invoking Lemma \ref{lemma3}, we will have
\[
\sum_{i=1}^K a_i\bigl(Y_{n}(\bu_i,\sigma_i)-\rE
Y_{n}(\bu_i,\sigma_i)\bigr)\Rightarrow N(0,\bga'\gS\bga)
\]
for any constants $a_{i}$, $i=1,\ldots,K$, where
\[
\bga=(a_1,\ldots,a_K)'
\]
and
\[
\gS_{ij}=\rE
Y(\bt_{i1},\bt_{i2},\sigma_i)Y(\bt_{j1},\bt_{j2},\sigma_j)=\vartheta(\bt
_{i1},\bt_{j2})\vartheta(\bt_{j1},\bt_{i2})
W(\sigma_i,\sigma_j)
\]
for the complex case and
\[
\gS_{ij}=
\bigl(\vartheta(\bt_{i1},\bt_{j2})\vartheta(\bt_{j1},\bt_{i2})
+\vartheta(\bt_{i1},\bt_{j1})\vartheta(\bt_{j2},\bt_{i2})\bigr)W(\sigma
_i,\sigma_j)
\]
for the real case.

To this end, we will verify the Liapounov condition and calculate
the asymptotic covariance matrix $\gS$ (see Lemma \ref{lemma3}) in
the next subsections.

\subsubsection{Verification of Liapounov's condition}

By (\ref{martina}), we have
%
%
\begin{eqnarray} \label{martinb1}\qquad
&&\sum_{i=1}^K
a_i\bigl(Y_{n}(\sigma_i)-\rE Y_{n}(\sigma_i)\bigr)\nonumber\\[-8pt]\\[-8pt]
&&\qquad=\sqrt{p}\sjln(\rE_j-\rE_{j-1}) \sum_{i=1}^K
(a_i\beta_j({\bx}^*(\bt_{i1})\bbA_j^{-1}(\sigma
_i)\bs_j\bs_j^*
\bbA^{-1}_j(\sigma_i)\bx(\bt_{i2}))).
\nonumber
\end{eqnarray}
The Liapounov condition with power index 4 follows by
verifying that
%
%
\begin{equation}\label{liap1}
p^2\sjln\rE\Biggl|\sum_{i=1}^K
a_i\beta_j{\bx}^*(\bt_{i1})\bbA_j^{-1}(\sigma_i)\bs_j\bs_j^*
\bbA^{-1}_j(\sigma_i)\bx(\bt_{i2})\Biggr|^4 \longrightarrow0.
\end{equation}
The limit (\ref{liap1}) holds if one can prove
that, for any ${\bx}_n,\by_n\in\Bbb C_1^p$,
%
%
\begin{equation}\label{liapa2}
p^2\sjln\rE| \beta_j{\bx}_n^*\bbA_j^{-1}(\sigma)\bs_j\bs_j^*
\bbA^{-1}_j(\sigma)\by_n|^4 \longrightarrow0 .
\end{equation}
To do this, applying Lemma 2.7 of \citet{BaiSil98}, for any $q\ge
2$, we get
%
%
\begin{equation}\label{ineq1}
\cases{\displaystyle\max_j\rE|\bs_j^*
\bbA_j^{-1}(\sigma)\by_n{\bx}_n^*\bbA_j^{-1}(\sigma)\bs_j|^q
=O(n^{-1-q/2}),\cr
\displaystyle\max_j\rE|\gma_j(\sigma)|^q=O(n^{-1-q/2}) \quad\mbox{and}
\cr\displaystyle\max_j\rE|\xi_j(\sigma)|^q=O(n^{-q/2}).}
\end{equation}
When $q>2$, the $O(\cdot)$ can be replaced by $o(\cdot)$ in the first
two inequalities. The assertion in (\ref{liapa2}) will then easily
follow from the estimations in (\ref{ineq1}) and the observation that
$|\beta_j(\sigma)|<1$.

\subsubsection{Simplification of $Y_n(\bu)-\rE Y_n(\bu)$}

For any ${\bx}_n,\by_n\in\Bbb C_1^p$, from (\ref{martina}), we have
%
%
\begin{eqnarray}\label{martina2}
&&{\bx}_n^* \bbA^{-1}(\sigma)\by_n-\rE{\bx}_n^*
\bbA^{-1}(\sigma)\by_n \nonumber\\
&&\qquad=
\sjln\bigl(\bar b \rE_j\gma_j+\rE_j(b_j-\bar b)\gma_j\\
&&\qquad\quad\hspace*{16.7pt}{}
+(\rE_j-\rE_{j-1})b_j(\sigma)\beta_j(\sigma)
\xi_j(\sigma)\bs_j^*
\bbA_j^{-1}(\sigma)\by_n{\bx}_n^*\bbA_j^{-1}(\sigma)\bs_j\bigr
).\nonumber
\end{eqnarray}
For the third term on the right-hand side of
(\ref{martina2}), applying (\ref{ineq1}), we have
\begin{eqnarray*}
&&\rE\Biggl|\sqrt{p}\sjln(\rE_j-\rE_{j-1})b_j(\sigma)\beta
_j(\sigma)
\xi_j(\sigma)\bs_j^*
\bbA_j^{-1}(\sigma)\by_n{\bx}_n^*\bbA_j^{-1}(\sigma)\bs_j\Biggr|^2\\
&&\qquad=p\sjln\rE|(\rE_j-\rE_{j-1})b_j(\sigma)\beta_j(\sigma)
\xi_j(\sigma)\bs_j^*
\bbA_j^{-1}(\sigma)\by_n{\bx}_n^*\bbA_j^{-1}(\sigma)\bs_j|^2\\
&&\qquad\le p\sjln(\rE| \xi_j(\sigma)|^4\rE|\bs_j^*
\bbA_j^{-1}(\sigma)\by_n{\bx}_n^*\bbA_j^{-1}(\sigma)\bs_j|^4)^{1/2}=o(n^{-1/2}).
\end{eqnarray*}
For the second term on the right-hand side of (\ref{martina2}),
we have
\begin{eqnarray*}
&&\rE\Biggl|\sqrt{p}\sjln\rE_j\bigl(b_j(\sigma)-\bar
b(\sigma)\bigr)\gma_j(\sigma)\Biggr|^2\\
&&\qquad\le p\sjln\bigl(\rE|b_j(\sigma)-\bar
b(\sigma)|^4\rE|\gma_j(\sigma)|^4\bigr)^{1/2}\\
&&\qquad= o(n^{-3/2})\cdot\Bigl(\max_j\rE|{\rtr\bbA_j^{-1}}(\sigma)-\rE
\rtr
\bbA^{-1}(\sigma)|^4\Bigr)^{1/2}\\
&&\qquad=o(n^{-1/2}),
\end{eqnarray*}
where the last step
follows from applying the martingale decomposition and the Burkholder
inequality and using the fact that
\[
|{\rtr\bbA_j^{-1}}(\sigma)-\rtr\bbA^{-1}(\sigma)|\le1/\sigma
\]
and
\[
\rE|{\rtr\bbA^{-1}}(\sigma)-\rE\rtr\bbA^{-1}(\sigma)|^4=O(n^2).
\]
Thus, we conclude that
%
%
\begin{equation} \label{simp2}
\sqrt{p}\bigl({\bx}_n^* \bbA^{-1}(\sigma
)\by_n-\rE
{\bx}_n^* \bbA^{-1}(\sigma)\by_n\bigr) =\sqrt{p}\sjln\bar b \rE
_j\gma_j+o_p(1) .
\end{equation}

\subsection{Asymptotic covariances}

To compute $\gS$, by the limiting property in (\ref{simp2}), we
only need to compute the limit
\[
\nu_{i,j} = \lim p\skln\bar
b(\sigma_{i})\bar
b(\sigma_{j})\rE_{k-1}\rE_k\gma_k(\bt_{i1},\bt_{i2},\sigma
_{i})\rE_k\gma_k(\bt_{j1},\bt_{j2},\sigma_{j}),
\]
in which, for any $i,k=1,\ldots, K$, we have
\begin{eqnarray*}
\gma_k(\bt_{i1},\bt_{i2},\sigma_i)
&=&\bs_k^*\bbA_k^{-1}(\sigma_i)\bx(\bt_{i2})\bx^*(\bt_{i1})
\bbA_k^{-1}(\sigma_i)\bs_k\\
&&{}-\frac{1}{n}\bx^*(\bt_{i1})\bbA_k^{-1}(\sigma_i)\bbA_k^{-1}(\sigma
_i)\bx(\bt_{i2}).
\end{eqnarray*}
By Lemma \ref{lema0}, we obtain $\bar b(\sigma)\to b(\sigma
)=1/(1+ym(\sigma))$,
Thus, we only need to calculate
%
%
\begin{equation} \label{variance1}
\nu_{i,j}=\lim_n p\skln b(\sigma_{i})
b(\sigma_{j})\rE_{k-1}\rE_k\gma_k(\bt_{i1},\bt_{i2},\sigma
_i)\rE_k\gma_k(\bt_{j1}.\bt_{j2},\sigma_{j}).
\end{equation}

For simplicity, we will use ${\bx},\by,\bu,\bbv,\sigma_1$
and $\sigma_2$ to
denote $\bx(\bt_{i1})$, $\bx(\bt_{i2})$, $\bx(\bt_{j1})$,
$\bx(\bt_{j2})$, $\sigma_{i}$ and $\sigma_{j}$. For
$\bbX=(X_1,\ldots,X_p)'$ of i.i.d. entries with mean 0 and variance~1,
and $\bbA=(\bbA_{ij})$ and $\bbB=(\bbB_{ij})$ to
be Hermitian matrices, the following equality holds:
\begin{eqnarray*}
&&
\rE(\bbX^* \bbA\bbX-\rtr\bbA)(\bbX^*
\bbB\bbX-\rtr\bbB)\\
&&\qquad=\rtr\bbA\bbB+|\rE X_{1}^2|^2\rtr\bbA\bbB^T+\sum
\bbA_{ii}\bbB_{ii}(\rE|X_1|^4-2-|\rE X_1^2|^2).
\end{eqnarray*}
Using this
equality, we get
%
%
\begin{eqnarray}\label{siglim1}
\nu&= &\lim\frac{pb(\sigma_{1}) b(\sigma
_2)}{n^2}\skln
\rtr\bigl(\rE_{k} \bbA_k^{-1}(\sigma_1)\by{\bx}^* \bbA
_k^{-1}(\sigma_1)\nonumber\\[-8pt]\\[-8pt]
&&\hspace*{13.9pt}\hspace*{86.6pt}{}\times
\rE_{k}\bbA_k^{-1}(\sigma_2)\bbv\bu^* \bbA
_k^{-1}(\sigma_2)\bigr)\nonumber
\end{eqnarray}
for the complex case and obtain
%
%
\begin{eqnarray}\label{siglim1real}
\nu&=& \lim
\frac{pb(\sigma_{1}) b(\sigma_2)}{n^2}\skln\rtr\bigl(\rE_{k}
\bbA_k^{-1}(\sigma_1)\by{\bx}^* \bbA_k^{-1}(\sigma
_1)\nonumber\\[-8pt]\\[-8pt]
&&\hspace*{101.1pt}{}\times \rE_{k}\bbA_k^{-1}(\sigma_2)(\bbv\bu^*+\bu
\bbv^*)
\bbA_k^{-1}(\sigma_2)\bigr)\nonumber
\end{eqnarray}
for the real case.

One could easily calculate the limit in (\ref{siglim1}) by applying the
method used in \citet{BaiMiaPan07} and by using the proof of their
equation (4.7). Therefore, we only need to calculate the limit of
%
%
\begin{eqnarray}
\label{seca6}
\hspace*{22pt}&&\frac{yb(\sigma_1)b(\sigma_2)}{n}\sum_{k=1}^n
\rtr\rE_k(\bbA_k^{-1}(\sigma_1)\by\bx^*\bbA_k^{-1}(\sigma_1))
\rE_k(\bbA_k^{-1}(\sigma_2)\bbv\bu^*\bbA_k^{-1}(\sigma_2))\nonumber\\
[-8pt]\\[-8pt]
\hspace*{22pt}&&\qquad=\frac{yb(\sigma_1)b(\sigma_2)}{n}\sum_{k=1}^n\rE_{k}
(\bx^*\bbA_k^{-1}(\sigma_1)\breve\bbA_k^{-1}(\sigma_2)\bbv
)(\bu^*\breve
\bbA_k^{-1}(\sigma_2) \bbA^{-1}_k(\sigma_1)\by),\nonumber
\end{eqnarray}
where
$\breve\bbA_k^{-1}(z_2)$ is similarly defined as
$\bbA_k^{-1}(\sigma_2)$ by using
$(\bs_1,\ldots,\bs_{k-1},\breve{\bs}_{k+1},\break\ldots,\breve
{\bs}_n)$
and by using the fact that
$\breve{\bs}_{k+1},\ldots,\breve{\bs}_n$ are i.i.d. copies of
$\bs_{k+1},\ldots,\bs_n$.

Following the arguments in \citet{BaiMiaPan07}, we only have
to replace their vectors $\bx_n$ and $\bx_n^*$
connected with $\bbA^{-1}_k(\sigma_1)$ by $\by$ and $\bx^*$
and replace
those connected with $\bbA^{-1}_k(\sigma_2)$ by
$\bbv$ and $\bu^*$, respectively. Going along with the same lines
from their (4.7) to (4.23), we obtain
%
%
\begin{eqnarray} \label{ff0}\quad
&&\rE_{k}
\bx^*\bbA_k^{-1}(\sigma_1)\breve\bbA_k^{-1}(\sigma_2)\bbv\bu^*\breve
\bbA_k^{-1}(\sigma_2)\bbA_k^{-1}(\sigma_1)\by\nonumber\\
&&\quad{}\times
\biggl[1-\frac{k-1}{n}\bar b(\sigma_1)\bar b(\sigma_2) \frac{1}{n}\rtr
T^{-1}(\sigma_2) T^{-1}(\sigma_1)\biggr]\nonumber\\[-8pt]\\[-8pt]
&&\qquad=\bx^*T^{-1}(\sigma_1)T^{-1}(\sigma_2)\bbv
\bu^*T^{-1}(\sigma_2)T^{-1}(\sigma_1)\by\nonumber\\
&&\qquad\quad{}\times
\biggl(1+\frac{k-1}{n}\bar b(\sigma_1)\bar b(\sigma_2)
\frac{1}{n}\rE_{k-1}\rtr(\bbA_k^{-1}(\sigma_1)\breve
\bbA_k^{-1}(\sigma_2))\biggr)+o_p(1)\nonumber
\end{eqnarray}
and
\[
\rE_{k} \rtr(\bbA_k^{-1}(\sigma_1) \breve\bbA_k^{-1}(\sigma_2))
=\frac{\rtr(T^{-1}(\sigma_1)T^{-1}(\sigma_2))+o_p(1)}{1-(({k-1})/{n^2})
b(\sigma_1)b(\sigma_2)\rtr(T^{-1}(\sigma_1)T^{-1}(\sigma_2))},
\]
where
\[
T(\sigma)=\biggl(\sigma+\frac{n-1}n b(\sigma)\biggr)\bbI.
\]
We then obtain
%
%
\begin{eqnarray}\label{ff1}
d(\sigma_1,\sigma_2):\!&=&\lim\bar b(\sigma
_1)\bar
b(\sigma_2) \frac{1}{n}\rtr(T^{-1}(\sigma_1)
T^{-1}(\sigma_2))\nonumber\\[-8pt]\\[-8pt]
&=&
\frac{y
b(\sigma_1)b(\sigma_2)}{(\sigma_1+b(\sigma_1))(\sigma_2+b(\sigma
_2))}\nonumber
\end{eqnarray}
and
%
%
\begin{eqnarray}\label{ff2}
\quad h(\sigma_1,\sigma_2):\!&=& b(\sigma_1) b(\sigma_2)
\bx^*T^{-1}(\sigma_1)T^{-1}(\sigma_2)\bbv
\bu^*T^{-1}(\sigma_2)T^{-1}(\sigma_1)\by\nonumber\\[-8pt]\\[-8pt]
&=&
\frac{\bx^*\bbv\bu^*\by b(\sigma_1)b(\sigma_2)}
{(\sigma_1+b(\sigma_1))^2(\sigma_2+b(\sigma_2))^2}.\nonumber
\end{eqnarray}
From
(\ref{ff1}) and (\ref{ff2}), we get
\begin{eqnarray*}
&&\mbox{The right-hand side of (\ref{seca6})}\\
&&\qquad\stackrel{\mathrm{a.s.}}\longrightarrow
yh(\sigma_1,\sigma_2)\biggl(\int_0^1\frac{1}{(1-td(\sigma
_1,z_2))}\,dt+\int_0^1\frac{td(\sigma_1,\sigma_2)}{(1-td(\sigma
_1,\sigma_2))^2}\,dt\biggr)
\\
&&\qquad=\frac{yh(\sigma_1,\sigma_2)}{1-d(\sigma_1,\sigma_2)}\\
&&\qquad=\frac{y\bx^*\bbv\bu^*\by b(\sigma_1)
b(\sigma_2)}{(\sigma_1+b(\sigma_1))(\sigma_2+b(\sigma_2))[(\sigma
_1+b(\sigma_1))(\sigma_2+b(\sigma_2))
-yb(\sigma_1) b(\sigma_2)]}.
\end{eqnarray*}
In addition, from (\ref{meq}), we establish
%
%
\begin{equation} \label{meqbb}
\frac{1}{\sigma+b(\sigma)}=m(\sigma) \quad\mbox{and}\quad
\frac{b(\sigma)}{\sigma+b(\sigma)}=1-\sigma m(\sigma).
\end{equation}
Applying
these identities, the limit of (\ref{seca6}) can be simplified to
\[
\bx^*\bbv\bu^*\by W(\sigma_1,\sigma_2),
\]
where
\[
W(\sigma_1,\sigma_2)=\frac{ym(\sigma_1)m(\sigma_2)}{1-y
(1-\sigma_1m(\sigma_1))(1-\sigma_2m(\sigma_2))}.\vadjust{\goodbreak}
\]
%
By symmetry, the limit of (\ref{siglim1real}) for the real case can
also be
simplified to
\[
(\bx^*\bbv\bu^*\by+\bx^*\bu\bbv^*\by)W(\sigma
_1,\sigma_2) .
\]
Therefore, for the complex case, the covariance function of the
process $Y(\bt_{i1}$, $\bt_{i2},\sigma)$ is
\[
\rE
Y(\bt_{i1},\bt_{i2},\sigma_1)Y(\bt_{j1},\bt_{j2},\sigma_2)=\vartheta(\bt
_{i1},\bt_{j2})\vartheta(\bt_{i2},\bt_{j1})
W(\sigma_1,\sigma_2),
\]
while, for the real case, it is
\begin{eqnarray*}
&&\rE
Y(\bt_{i1},\bt_{i2},\sigma_1)Y(\bt_{j1},\bt_{j2},\sigma_2)\\
&&\qquad=\bigl(\vartheta(\bt_{i1},\bt_{j2})\vartheta(\bt_{j1},\bt_{i2})
+\vartheta(\bt_{i1},\bt_{j1})
\vartheta(\bt_{j2},\bt_{i2})\bigr)W(\sigma_1,\sigma_2).
\end{eqnarray*}


\subsection{Tightness}\label{sectight}
\begin{theorem}\label{them4}
Under the conditions in Theorem \ref{thm1}, the sequence of
$Y_n(\bu,\sigma)-\rE(Y_n(\bu,\sigma)$ is tight.
\end{theorem}

For ease reference on the tightness, we quote a proposition
from page 267 of \citet{Loe78} as follows.
\begin{proposition}[(Tightness criterion)] \label{prop1} The sequence
$\{\rP_n\}$
of probability measure is tight if and only if:
\begin{longlist}[(ii)]
\item[(i)]
\begin{eqnarray*}
\\[-40pt]
\sup_n \rP_n\bigl(x\dvtx|x(0)|>c\bigr) &\longrightarrow&0
\qquad\mbox{as }c \to\infty
\end{eqnarray*}
and, for every $\ep>0$, as $\gd\to0$, we have\vspace*{6pt}
\item[(ii)]
\begin{eqnarray*}
\\[-40pt]
\rP_n\bigl(\omega_x(\gd)>\ep\bigr) &\longrightarrow&0,
\end{eqnarray*}
where $\gd$-oscillation is defined by
\[
\omega_x(\gd)=\sup_{|\bt-\bs|<\gd}|x(\bt)-x(\bs)|.
\]
\end{longlist}
\end{proposition}

To complete the proof of the tightness for Theorem \ref{them4}, we
note that condition (i) in Proposition \ref{prop1} is a consequence
of finite-dimensional convergence which has been proved in the
previous section. To demonstrate condition (ii)
in Proposition~\ref{prop1}, we will use the two lemmas given below.
Therefore, to complete the proof of Theorem \ref{them4}, by
Proposition \ref{prop1} and Lemma \ref{lemtight1}, it is sufficient
to verify that
%
%
\begin{equation}\label{finite}
\sup_{\bu_1,\bu_2\in T\times
T}\rE\biggl|\frac{Y_n(\bu_1)-Y_n(\bu_2)}{\|(\bu_1,\sigma
_1)-(\bu_2,\sigma_2)\|}\biggr|^{4m+2}<\infty.
\end{equation}
This inequality will be proved in Lemma \ref{lemtight2} stated below.
\begin{lemma}\label{lemtight1} Suppose that
$X_n(\bt)$ is a sequence of stochastic processes, defined on an
$m$-dimensional time domain $T$,
whose paths
are continuous and Lipschitz; that is, there is a random variable
$R=R_n$ such that
\[
|X_n(\bt)-X_n(\bs)|\le R\|\bt-\bs\|.
\]
If there is an $\alpha>m$ such that
%
%
\begin{equation}
\label{tight3}
\rE|R|^\ga<\infty,
\end{equation}
then, for any fixed $\ep>0$, we have
%
%
\begin{equation}\label{eposci}
\lim_{\gd\downarrow0}\rP_n\bigl(\omega_x(\gd)>\ep\bigr)= 0.
\end{equation}
\end{lemma}
\begin{pf}
Without loss of generality, we assume that $T=[0,M]^m$.
First, for any given $\ep>0$ and $\gd>0$, we choose an integer $K$
such that $MK^{-1}<\gd$ and $2^\ga K^{m-\ga}<1/2$.
For each $\ell=1,2,\ldots,$ we define
\[
t_{i}(j,\ell)=\frac{jM}{K^\ell},\qquad j=1,\ldots,K^\ell.
\]
Denoting by $\bt(\bbj,\ell)$, $\bbj=(j_1,\ldots,j_m)$, the vector
whose $i$th entry is $t_{i}(j_i,\ell)$. Then, we have
\begin{eqnarray*}
\hspace*{-4pt}&& \rP_n\bigl(\omega_x(\gd)\ge\ep\bigr)\\
\hspace*{-4pt}&&\qquad\le
2\rP\Bigl(\sup_{\bbj,1} \sup_{|\bt-\bt(\bbj,1)|\le
2M\sqrt{m}/K}|X_{n}(\bt)-X_n(\bt(\bbj,1))|\ge\ep/2\Bigr)\\
\hspace*{-4pt}&&\qquad\le\sum_{\ell=1}^L \sum_{(\bbj,\ell+1)}2\rP\bigl(
\bigl|X_{n}(\bt(\bbj^*,\ell))-X_n\bigl(\bt(\bbj,\ell+1)\bigr)\bigr|\ge
2^{-\ell-1}\ep\bigr)\\
\hspace*{-4pt}&&\qquad\quad{}+2\rP\Bigl(\sup_{\bt(\bbj,L+1)}\sup_{\|\bt
-\bt(\bbj
,L+1)\|\le
2\sqrt{m}M/K^{L}}
\bigl|X_{n}(\bt)-X_n\bigl(\bt(\bbj,L+1)\bigr)\bigr|\ge
2^{-L-2}\ep\Bigr) \\
\hspace*{-4pt}&&\qquad\le\sum_{\ell=1}^\infty2(K^\ell/M)^m\biggl(\frac
{2\sqrt
{m}M}{\ep
2^{-\ell-1}K^{\ell}}\biggr)^{\ga}\rE|R|^\ga\\
\hspace*{-4pt}&&\qquad=2^{2+3\ga}\ep^{-\ga}m^{\ga/2}(M/K)^{\ga-m}\rE
|R|^\ga\\
\hspace*{-4pt}&&\qquad=2^{2+3\ga}\ep^{-\ga}m^{\ga/2}\gd^{\ga-m}\rE
|R|^\ga\to0
\qquad\mbox{as } \gd\to0,
\end{eqnarray*}
where the summation $\sum_{(\bbj,\ell+1)}$ runs
over all possibilities of $j_i\le K^{\ell+1}$, and
$\bt(\bbj^*,\ell)$ is the $\bt(\bbj,\ell)$ vector closest to
$\bt(\bbj,\ell+1)$. Here, to prove the first
inequality, one only needs to choose $\bt(\bbj,1)$ as the center of
the first layer hypercube in which $\frac12(\bt+\bs)$ lies. The
proof of the second inequality could be easily obtained by applying a simple
induction. In the proof of the third inequality, the first term follows
by the
Chebyshev inequality and the fact that
\[
\bigl|X_{n}(\bt(\bbj^*,\ell))-X_n\bigl(\bt(\bbj,\ell+1)\bigr)\bigr|\le
R\|\bt(\bbj^*,\ell)-\bt(\bbj,\ell+1)\|\le R\sqrt
{m}M/K^{\ell}.
\]
At the same time, the second term tends to 0 for all fixed $n$ when
$L\to\infty$ because
\begin{eqnarray*}
&&\rP\Bigl(\sup_{\bt(\bbj,L+1)}\sup_{\|\bt-\bt(\bbj
,L+1)\|\le
2M\sqrt{m}K^{-L-2}} \bigl|X_{n}(\bt)-X_n\bigl(\bt(\bbj
,L+1)\bigr)\bigr|\ge
2^{-L-2}\ep\Bigr)
\\
&&\qquad\le\rP\bigl(|R|\ge(K/2)^{L+2}\ep/2M\sqrt{m}\bigr)\to0.
\end{eqnarray*}
Thus,
the proof of the lemma is complete.
\end{pf}
\begin{lemma}\label{lemtight2}
Under the conditions of Theorem \ref{thm1}, the property in
(\ref{finite}) holds for any $m$.
\end{lemma}
\begin{pf}
For simplicity, we only prove the
lemma for a general $m$ instead of $4m+2$. For a constant $L$, we
have
\begin{eqnarray*}
&&\rE\biggl|\frac{Y_n(\bu_1,\sigma_1)-Y_n(\bu_2,\sigma_2)-\rE
(Y_n(\bu_1,\sigma_1)-Y_n(\bu_2,\sigma_2))}
{\|\bu_1-\bu_2\|+|\sigma_1-\sigma_2|}\biggr|^m
\\
&&\qquad\asymp p^{m/2}\rE\biggl|
\frac{\bx_n(\bt_1)^*\bbA^{-1}(\sigma_1)\bx_n(\bt_2)-\rE\bx_n(\bt_1)^*\bbA
^{-1}(\sigma_1)\bx_n(\bt_2)}
{\|\bt_1-\bt_3\|+\|\bt_2-\bt_4\|+|\sigma_1-\sigma
_2|}\\
&&\phantom{\asymp p^{m/2}\rE\biggl|}
\qquad{}-\frac{(\bx_n(\bt_3)^*\bbA^{-1}(\sigma_2)\bx_n(\bt_4)-\rE\bx_n(\bt
_3)^*\bbA^{-1}(\sigma_2)
\bx_n(\bt_4))}{\|\bt_1-\bt_3\|+\|\bt_2-\bt_4\|
+|\sigma_1-\sigma_2|}\biggr|^m
\\
&&\qquad\le Ln^{m/2}\bigl\{\rE\bigl|\bigl(\bigl(\bx_n(\bt_1)-\bx_n(\bt_3)\bigr)^*\bbA^{-1}(\sigma
_1)\bx_n(\bt_2)\\
&&\qquad\quad\hspace*{46.7pt}{}-\rE\bigl(\bx_n(\bt_1)-\bx_n(\bt_3)\bigr)^*\bbA^{-1}(\sigma
_1)\bx_n(\bt_2)\bigr)({\|\bt_1-\bt_3\|})^{-1}\bigr|^m
\\
&&\phantom{Ln^{m/2}\bigl\{}
\qquad\quad{}+\rE\bigl|\bigl(\bx_n^*(\bt_3)\bbA^{-1}(\sigma_1)\bigl(\bx_n(\bt_2)-\bx_n(\bt
_4)\bigr)\\
&&\qquad\quad\hspace*{59.5pt}{}-\rE\bx_n^*(\bt_3)
\bbA^{-1}(\sigma_1)\bigl(\bx_n(\bt_2)-\bx_n(\bt_4)\bigr)\bigr)(\|
\bt_2-\bt_4\|)^{-1}\bigr|^m\\
&&\phantom{Ln^{m/2}\bigl\{}
\qquad\quad{}
+\rE|\bx_n^*(\bt_3)\bbA^{-1}(\sigma_1)\bbA
^{-1}(\sigma_2)\bx_n(\bt_4)\\
&&\qquad\quad\hspace*{133pt}{}-\rE\bx_n^*(\bt_2)
\bbA^{-1}(\sigma_1)\bbA^{-1}(\sigma_2)\bx_n(\bt_4)
|^m\bigr\},
\end{eqnarray*}
where $a\asymp b$ means $a$ and $b$ have the same order, that is,
there exists a positive constant $K$ such that $K^{-1}b<a<Kb$.

We note that $\|\bx_n(\bt_1)-\bx_n(\bt_3)\|/\|\bt_1-\bt_3\|\le
1$ or bounded for the general case. By applying the martingale
decomposition in (\ref{martina}), the Burkholder inequality and the
estimates in (\ref{ineq1}), we have
\begin{eqnarray*}
&&n^{m/2}\rE\biggl|\frac{(\bx_n(\bt_1)-\bx_n(\bt_3))^*\bbA^{-1}(\sigma_1)\bx
_n(\bt_2)
-\rE(\bx_n(\bt_1)-\bx_n(\bt_3))^*\bbA^{-1}(\sigma
_1)\bx_n(\bt_2)}{\|\bt_1-\bt_3\|}\biggr|^m\\
&&\qquad=O(1).
\end{eqnarray*}
Similarly, we obtain
\begin{eqnarray*}
&&n^{m/2}\rE\biggl|\frac{\bx_n^*(\bt_3)\bbA^{-1}(\sigma
_1)(\bx_n(\bt_2)-\bx_n(\bt_4))-\rE\bx_n^*(\bt_3)
\bbA^{-1}(\sigma_1)(\bx_n(\bt_2)-\bx_n(\bt_4))}{\|
\bt_2-\bt_4\|}\biggr|^m\\
&&\qquad=O(1).
\end{eqnarray*}
Using the martingale decomposition and the Burkholder inequality, we get
\begin{eqnarray*}
&&n^{m/2}\rE|\bx_n^*(\bt_3)\bbA^{-1}(\sigma_1)\bbA
^{-1}(\sigma_2)\bx_n(\bt_4)-\rE\bx_n^*(\bt_2)
\bbA^{-1}(\sigma_1)\bbA^{-1}(\sigma_2)\bx_n(\bt_4)
|^m\\
&&\qquad\le Ln^{m/2}\Biggl[\skln
\rE|\bx_n^*(\bt_3)[\bbA^{-1}(\sigma_1)\bbA^{-1}(\sigma
_2)-\bbA_{k}^{-1}(\sigma_1)\bbA_{k}^{-1}(\sigma_2)]\bx_n(\bt_4)|^m\\
&&\qquad\quad\hspace*{32.3pt}{}+\rE\Biggl(\skln
\rE_{k-1}|\bx_n^*(\bt_3)[\bbA^{-1}(\sigma_1)\bbA
^{-1}(\sigma_2)\\
&&\qquad\quad\hspace*{129pt}{}-\bbA_{k}^{-1}(\sigma_1)\bbA_{k}^{-1}(\sigma
_2)]\bx_n(\bt_4)|^2
\Biggr)^{m/2}\Biggr]\\
&&\qquad= O(1),
\end{eqnarray*}
which follows from applying the following
decomposition:
\begin{eqnarray*}
\hspace*{-5pt}&&\bbA^{-1}(\sigma_1)\bbA^{-1}(\sigma_2)-\bbA_{k}^{-1}(\sigma
_1)\bbA_{k}^{-1}(\sigma_2)\\
\hspace*{-5pt}&&\qquad=\beta(\sigma_1)\bbA^{-1}_k(\sigma_1)\bs_k\bs_k^*
\bbA^{-1}_k(\sigma_1)\bbA_k^{-1}(\sigma_2)+
\beta(\sigma_2)\bbA^{-1}_k(\sigma_1)\bbA_k^{-1}(\sigma_2)\bs_k\bs_k^*
\bbA^{-1}_k(\sigma_2)\\
\hspace*{-5pt}&&\qquad\quad{}+\beta(\sigma_k)\beta(\sigma_2)\bbA^{-1}_k(\sigma_1)\bs_k\bs_k^*
\bbA^{-1}_k(\sigma_1)\bbA_k^{-1}(\sigma_2)\bs_k\bs_k^*
\bbA^{-1}_k(\sigma_2)
\end{eqnarray*}
and thereafter employing
the results in (\ref{ineq1}). Thus, condition (\ref{finite}) is
verified.
\end{pf}

\section{\texorpdfstring{Proof of Corollary \protect\ref{coro2}}{Proof of Corollary 2}}
\label{sec:3}

Applying the quadratic equation (\ref{meq}), we have
%
%
\begin{equation} \label{meq2}
\sigma=\frac1m-\frac1{1+ym}.
\end{equation}
Making a difference of
$\sigma_1$ and $\sigma_2$, we obtain
\[
\sigma_1-\sigma_2=\frac{m(\sigma_2)-m(\sigma_1)}{m(\sigma
_1)m(\sigma_2)}-\frac{y(m(\sigma_2)-m(\sigma_1))}{(1+ym(\sigma
_1))(1+ym(\sigma_2))}.
\]
We also establish
%
%
\begin{equation}\label{mdiff}
\frac{m(\sigma_2)-m(\sigma_1)}{\sigma_1-\sigma_2}=
\frac{m(\sigma_1)m(\sigma_2)(1+ym(\sigma_1))(1+ym(\sigma
_2))}{(1+ym(\sigma_1))(1+ym(\sigma_2))-ym(\sigma_1)m(\sigma_2)}.
\end{equation}
Finally, we conclude that
\begin{eqnarray*}
&&\frac{m(\sigma_2)-m(\sigma_1)}{\sigma_1-\sigma_2}-m(\sigma
_1)m(\sigma_2) \\
&&\qquad=
\frac{ym^2(\sigma_1)m^2(\sigma_2)}{(1+ym(\sigma_1))(1+ym(\sigma
_2))-ym(\sigma_1)m(\sigma_2)}=W(\sigma_1,\sigma_2)
\end{eqnarray*}
by noticing that $1+ym(\sigma)=m(\sigma)/(1-\sigma m(\sigma))$,
an easy
consequence of (\ref{meq2}).

Furthermore, one could easily show that the left-hand side of the
above equation is
\[
\int_a^b\frac{dF_y(x)}{(x+\sigma_1)(x+\sigma_2)}
-\int_a^b\frac{dF_y(x)}{x+\sigma_1}\int_a^b\frac{dF_y(x)}{x+\sigma_2}.
\]
By using the unique extension of analytic functions, we have
\[
W(z_1,z_2)=\int_a^b\frac{dF_y(x)}{(x-z_1)(x-z_2)}
-\int_a^b\frac{dF_y(x)}{x-z_1}\int_a^b\frac{dF_y(x)}{x-z_2} .
\]
Substituting this into Corollary \ref{coro1}, we complete the proof of
Corollary \ref{coro2}.

\begin{appendix}\label{app}
\section*{Appendix}

\begin{lemma}[{[Theorem 35.12 of Billingsley (\citeyear{Bil68})]}]\label
{lemma3} Suppose
that, for each n, $X_{n,1},X_{n,2},\ldots, X_{n,r_n}$ is a real
martingale
difference sequence with respect to the increasing $\sigma$-field
$\{{\mathcal{F}}_{n,j}\}$ having second moments. If, as $n\rightarrow
\infty$,
\begin{eqnarray*}
\mbox{\textup{(i)}\hspace*{6.78pt}}\quad
\sum_{j=1}^{r_n}\rE(X^2_{n,j}|{\mathcal
{F}}_{n,j-1})&\stackrel{\mathit{i.p.}}{\longrightarrow}&\sigma^2
\quad\mbox{and}\\
\mbox{\textup{(ii)}}\quad
\sum_{j=1}^{r_n}\rE\bigl(X^2_{n,j}I_{(|Y_{n,j}|\geq\varepsilon
)}\bigr)&\longrightarrow&0,
\end{eqnarray*}
where $\sigma^2$ is a positive constant and $\varepsilon$ is an
arbitrary positive number, then
\[
\sum_{j=1}^{r_n}X_{n,j}\stackrel{\mathscr D}\longrightarrow
N(0, \sigma^2).
\]
\end{lemma}

In what follows, $\bs_j$, $\bbA^{-1}$ and $\bbA_j^{-1}$ are
defined in Section \ref{sec:2} and $\bbM_j$ and $\bbM$ refer to any
pair of matrices which are independent of $\bs_j$.
\begin{lemma}\label{lema0} Under the conditions of Theorem
\ref{thm1}, for any matrix $M_j$ bounded in norm and independent of
$\bs_j$, we have
%
\setcounter{equation}{0}
\begin{equation} \label{a00}
\max_j\biggl|\fn(\bs_j^*M_j\bs_j-\rtr M_j)\biggr|\stackrel{\mathit{a.s.}}
\longrightarrow0.
\end{equation}
\end{lemma}

The proof of this lemma could be easily obtained by applying the
truncation technique and invoking Lemma 2.7 of \citet{BaiSil98}.
\begin{lemma}\label{lema01} Under the conditions of Theorem
\ref{thm1}, for any $\bx_n,\by_n\in\Bbb C_1^p$,
%
%
\begin{equation} \label{a01}
\sup_j|\bx_n^*\bbA^{-1}M\by_n-\bx_n^*\bbA_j^{-1}M\by_n|\stackrel{\mathit{a.s.}}
\longrightarrow0.
\end{equation}
Similarly, for any matrix $M$
with bounded norm and independent of $\bs_i$, we have
%
%
\begin{equation}\label{a02}
\max_{i,j}|\rE_j\bx_n^*\bbA_j^{-1}(\sigma)M\by_n-\rE_j\bx_n^*\bbA
_{ij}^{-1}(\sigma)M\by_n|
\stackrel{\mathit{a.s.}}\longrightarrow0.
\end{equation}
\end{lemma}
\begin{pf}
Using
%
%
\begin{equation}\label{ajinv}
\bbA^{-1}(\sigma)=\bbA^{-1}_j(\sigma)-\bbA^{-1}_j(\sigma)\bs_j\bs
_j^*\bbA^{-1}_j(\sigma)\beta_j(\sigma),
\end{equation}
we obtain
\begin{eqnarray*}\sup_j|\bx_n^*\bbA^{-1}M\by_n-\bx_n^*\bbA
_j^{-1}M\by_n|
&\le&
\sup_j\fn|\bx_n^*\bbA_j^{-1}\bs_j\bs_j^*\bbA
_j^{-1}M\by_n| \\
&=& \sup_j\fn|\bx_n^*\bbA_j^{-1}\bbA_j^{-1}M\by_n|+o(1),
\end{eqnarray*}
which, in turn, implies (\ref{a01}). Here, we adopt (\ref{a00}) in the
last step above. The conclusion (\ref{a02}) can be proved in a similar way.
\end{pf}
\begin{lemma}\label{lema1} Under the conditions of Theorem
\ref{thm1}, for any $\bx_n,\by_n\in\Bbb C_1^p$, we have
%
%
\begin{equation} \label{a11}
\bx_n^*\bbA^{-1}(\sigma)\by_n-\bx^*_n\by_n
m(\sigma)\stackrel{\mathit{a.s.}}
\longrightarrow0
\end{equation}
and
%
%
\begin{equation}\label{a12}
\max_j|\rE_j\bx_n^*\bbA_j^{-1}(\sigma)\by_n-\bx_n^*\by_n
m(\sigma)|\stackrel{\mathit{a.s.}} \longrightarrow0.
\end{equation}
\end{lemma}
\begin{pf}
By using the formula $\bbA=\sjln
\bs_j\bs_j^*+\sigma\bbI$ and multiplying $\bx_n^*$ from
the left- and multiplying $\bbA^{-1}\by_n$ from the right-hand side of the
equation, we obtain
\[
\bx_n^*\bbA^{-1}(\sigma)\by_n=\sigma^{-1}\bx_n^*\by_n-\frac1{n\sigma
}\sjln
\bx_n^*\bs_j\bs_j^*\bbA_j^{-1}(\sigma)\beta_j(\sigma
)\by_n.
\]
As $\beta_j(\sigma)\stackrel{\mathrm{a.s.}} \longrightarrow b=\frac
1{1+ym(\sigma)}$
uniformly in $j$, we apply Lemmas \ref{lema0} and \ref{lema01} and
obtain
\[
\bx_n^*\bbA^{-1}(\sigma)\by_n=\sigma^{-1}\bx_n^*\by_n-\sigma^{-1}
\bx_n^*\bbA^{-1}(\sigma)\by_nb(\sigma)+o(1).
\]
This, in turn, implies that
\[
\bx_n^*\bbA^{-1}(\sigma)\by_n=\frac{\bx_n^*\by_n+o(1)}{\sigma
+b(\sigma)}.
\]
The conclusion in (\ref{a11}) could then follow from the fact that
\[
m(\sigma)=\frac{1}{\sigma+b(\sigma)},
\]
whereas the conclusion in (\ref{a12}) can be proved by employing the
same method.
\end{pf}
\begin{lemma}\label{lem33} Under the conditions of Theorem
\ref{thm1}, for any $\bx_n,\by_n\in\Bbb C_1^p$, we have
\[
\sqrt{n}\bigl(\bx_n^*\rE(\bbS_n+\sigma
I)^{-1}\by_n-\bx_n^*\by_n m_n(\sigma)\bigr)\longrightarrow0.
\]
\end{lemma}
\begin{pf}
When $\by_n=\bx_n$, Lemma \ref{lem33} in our paper reduces
to the conclusion (5.5) $\to0$ as shown in \citet{BaiMiaPan07}.
To complete the proof, one could simply keep
$\bx_n^*$ unchanged and substitute $\bx_n$ by
$\by_n=(\bx_n^*\by_n)\bx_n+\bbz_n$ in the proof of the above
conclusion. Thereafter, the proof of this lemma follows.
\end{pf}
\end{appendix}

\section*{Acknowledgments}
The authors are grateful to the Editor, Professor Andrew Barbour,
Associate Editor, Professor Rick Durrett and an anonymous referee for
their substantive comments and suggestions that have significantly
improved the manuscript. We would also like to show our appreciation to
Ms. Sarah A. Burke and Ms. Ke Xin Tan for their assistance in editing
our paper.

%

%
\printaddresses

\end{document}